# A novel multi-stage concurrent topology optimization for variable-stiffness structures

Yaya Zhang[a] [*], Hu Wang[a,b] [†], Jichao Yin[a,b], Shuhao Li[a], Mengzhu Yang[c]

[a] *State Key Laboratory of Advanced Design and Manufacturing for Vehicle Body, Hunan University, Changsha 410082, People's Republic of China;*

[b] *Beijing Institute of Technology Shenzhen Automotive Research Institute, Shenzhen 518000, People's Republic of China*

[c] *Department of Aerospace Engineering, University of Bristol Bristol, BS8 1TR United Kingdom;*

[*]Hu Wang: wanghu@hnu.edu.cn

## Highlights

- A new scheme of concurrent optimization of topology and fibre orientation is proposed.

- The SBPTO is suggested to ensure the fibre angle optimization convergence.

## Abstract

The concurrent optimization of topology and fibre orientation is a promising approach to pursue higher strength and lighter weight of variable-stiffness structure. This study proposes a novel discrete-continuous scheme for the concurrent optimization. Considering the global convergence, Discrete Material Optimization (DMO) is firstly utilized to select a dominant angle for each element from several predefined candidate angles. However, it is still difficult to obtain excellent fibre angle convergence due to difficulty in selection of candidate materials for some elements. Therefore, the

[*] First author. E-mail address: yyzh200@hnu.edu.cn (Y. Zhang)

[†] Corresponding author. E-mail address: wanghu@hnu.edu.cn (H. Wang)

Sequential Binary-Phase Topology Optimization is proposed to guarantee the uniqueness of the element to candidate angle mapping. Moreover, to obtain better mechanical properties, Continuous Fibre Angle Optimization (CFAO) with spatial filter is introduced to optimize fibre continuity. Several classic numerical examples are used to verify the excellent performance of the proposed method in terms of fibre angle convergence and stable optimization ability.

# 1. Introduction

Structure lightweight design plays an important role in aerospace [1, 2], shipping [3] and automotive tools[4]. Generally, structure lightweight design includes two main configurations[5]. Fibre-reinforced composites have been widely used in these fields due to the excellent mechanical properties of high strength, stiffness and light weight [6, 7] . The fibre reinforced composites can be divided into Constant Stiffness (CS) composites and variable stiffness (VS) composites. CS composites, in which the fibre angle of each layer is constant, are typically optimized by the laminate stacking sequence, lamination parameters and fibre angle of each layer [8-10]. Compared with the CS, VS composites with variable fibre angles allows structures to distribute the loads more efficiently [11]. Its optimization objects typically include lamination parameters, fibre path and local fibre angles. Automated fibre placement(AFP) [12], filament winding [13] and automated tape laying (ATL) [14] enable fabrication of curvilinear fibres and that makes it possible to fabricate the VS composites. The VS optimization to improve load carrying capacity as well as light structure weight has gradually attracted the attention of researchers [15].

Topology Optimization (TO) has excellent performance in seeking new structural configurations while reducing weight. TO guides engineers to place materials in the prescribed design domain for the excellent structural performance [16]. Compared to other optimization methods, TO provides a greater degree of freedom in the design space and is an effective design method for structural light-weighting. Topology optimization considering fibre angle can significantly improve the structural stiffness [17, 18], natural frequency [19] and maximum buckling load [20]. Therefore, the combination of structural topology and fibre angle orientation optimization is a very promising approach for structure lightweight design in fibre reinforced composites.

In fibre orientation optimization, the continuity of fibre orientation plays a crucial role in ensuring material properties such as strength or stiffness because fibre discontinuity leads to stress concentration [21]. The continuous fibre angle optimization faces a challenge that the rotated stiffness tensor is composed of multi-valued functions, such as trigonometric function. The periodicity of trigonometric functions makes the optimization problem a highly non-convex problem with multiple local optimal solutions in given space. To avoid local convergence, non-gradient optimization--Heuristic Optimization Algorithms (HOA) had been widely employed, such as Genetic Algorithms (GA) [22, 23], Artificial Immune Algorithms (AIA) [24], Ant Colonies (AC) [25] and Particle Swarm Optimization (PSO) [26]. Keller [27] adopted the Evolutionary Algorithm with first-order search and the niching strategy to optimize the orientation angles. Other studies [28-30] used GA to optimize the orientation angles, number of plies, and stacking sequence. In these studies, the ply angles were assumed to be constant. However, Sigmund [31] compared Non-Gradient Topology Optimization (NGTO) [32] with Gradient-based Topology Optimization (GTO) from the perspective of global search ability and computational efficiency. The results shows that NGTO commonly could not find global optima solutions based on the coarse meshes and have high computational cost. Furthermore, NGTO is difficult to handle the problems of high-aspect ratio, complex geometry design domains, complex physics situations and accuracy demands. They proved that GTO performs much better in these areas than NGTO.

For the variable stiffness design problem, several gradient-based methods have been proposed to optimize the fibre angle orientation. Based on discrete Kirchhoff theory, Mota Soares et. al [33] carried out sensitivity analysis of fibre orientation as well as the thickness for the optimization. Bruyneel et. al [34] adopt the approximation concepts approach and dual algorithm. They performed sensitivity analysis analytically and used sequential convex programming--Method of Moving Asymptotes (MMA) [35], to optimize the fibre orientation and the thickness efficiently. This method is also called Continuous Fibre Angle Optimization (CFAO). CFAO takes the angle directly as a design variable, which varies in the whole angle orientation range. It tends to fall into local optima and the optimization result is highly dependent on the initial design values [36, 37].

To solve these problems, Stegmann and Lund [38] suggested a method, Discrete Material Optimization (DMO), whose equivalent constitutive matrix is expressed as a weighted sum of candidate constitutive matrixes from prescribed discrete angles. DMO based on gradient and penalty coefficient is also adopted to force each element to correspond to only one candidate angle. Experimental results show that the DMO has good global optimization ability and is not sensitive to the initial design values. Subsequently, Lund [20] studied the problem of maxing buckling load factor of multi-material composite structures by DMO. Bruyneel [39] proposed a new variant of DMO, the Shape Functions with Penalization scheme (SFP), with fewer design variables required. Yin et al. [40] introduced the peak function into material interpolation for optimal selection of different isotropic materials which advantage is that a variety of materials can be introduced into optimization without increasing design variables. However, artificially high stiffness materials may occur during the process, which may lead to local minimum. Gao et al. [41] proposed the Bi-value Coding Parameterization (BCP) which can reduce the number of design variables significantly. Wu et al. [6] combined DMO with commercial software ABAQUS to optimize the ply angles of the laminate vehicle door.

Compared with CFAO, DMO has better global optimization ability and low sensitivity to the initial design values. But because of fibre discontinuities, DMO has poor manufacturability and will easily lead to stress concentration. To overcome the shortcomings of CAFO and DMO, Kiyono et al. [42] proposed a novel scheme of discrete-continuous fibre angle optimization. This study used the normal distribution functions as weighting functions in DMO to choose one discrete angle from any number of discrete candidate angles for each element. And then a spatial filter was utilized to ensure the continuity of the fibres, where the filter radius would influence the level of fibre continuity. Luo et al. [43] proposed a concurrent optimization framework for topology optimization and discrete-continuous fibre angle optimization. In the study, the fibre orientation interval is divided into several average subintervals. The optimization problem [10] becomes to select an interval among several discrete subintervals and perform continuous angle optimization within the interval. The DMO and CFAO have been carried out in sequence in each iteration step. Ding. and Xu. introduced the normal distribution functions as the weight functions to the concurrent optimization. Researchers have proposed some good ideas for the concurrent optimization, but there are still several problems to be addressed. Firstly, despite the use of penalty strategies, the elements convergence is still not guaranteed and the convergence process is relatively slow. Secondly, the number of subintervals as well as the number of discrete design variables affects the optimization result.

Considering the advantages of DMO and CFAO, this paper seeks a new scheme for the concurrent optimization of structural topology and fibre angle orientation. Due to good global optimization ability, DMO has been utilized to select an optimized angle for each element from several predefined candidate angles firstly. However, some "fuzzy regions" make it difficult to select the best angle for each element. To avoid "fuzzy regions", a Sequential Binary-phase Topology Optimization (SBPTO) is suggested to optimize further based on the optimized structure of the DMO. The SBPTO is a topology optimization method for solving multiphase problem where one phase represents a candidate angle. The essential of the SBPTO is to decompose the multiphase topology optimization problem into a series of two-phase topology optimization problem in a sequential manner. This method can be generally and easily implemented. The integration of DMO and SBPTO not only preserves the optimization ability of DMO, but also solves the problem of fibre convergence for the interpolation penalty model and improves the stability of the solution. Finally, considering mechanical properties, the continuity of the fibre path is optimized by using CFAO with spatial filter. The process of DMO-SBPTO-CFAO is abbreviated as DSCO (See **Fig. 1**).

The remainder of this paper is organized as follows: section 2 introduces the discrete-continuously fibre optimization method. In section 3, some numerical examples are shown. Finally, the findings are concluded in section 4. The code and data accompanying this study will be made publicly available at https://github.com/HnuAiSimOpt/DSCO after the paper is officially published.

## 2. Discrete-continuously fibre optimization method

The framework of the suggested concurrent optimization DSCO is illustrated in Fig. 1. In the proposed framework, DMO, SBPTO and CFAO are integrated. Considering the global convergence of DMO [37], the DMO is employed to find the "near optimum" solution efficiently. However, there might be some unavoidable "fuzzy regions" which make it difficult to select the best angle. According to Fig. 1, such "fuzzy regions" are usually located at the junction of different angles or near where constraints and forces are applied. The elements located in such regions are called unconvergent elements. In order to remove the unconvergent elements, the SBPTO is introduced to facilitate selecting the best angle from the candidate angles. The SBPTO decomposes the multi-phase topology optimization method into a series of binary phase topology optimization where each phase presents a candidate fibre angle or void-material. Based on the design space of DMO, the design variables of the elements that satisfy the conditions which would be explained in section 2.2 will be frozen. As the design space is reduced, the contribution of unconvergent elements becomes more prominent. Binary phase topology optimization makes it easier for the difficult-to-choose elements to converge. The CFAO is sensitive to the initial value and easily fall into local optimization due to the high nonconvex character with multiple local optimization [15, 29, 39]. The two-step discrete fibre optimization provides reliable initial design values for CFAO to reduce the risk of falling into local optimality. The flow chart of the concurrent optimization is as Fig. 2.

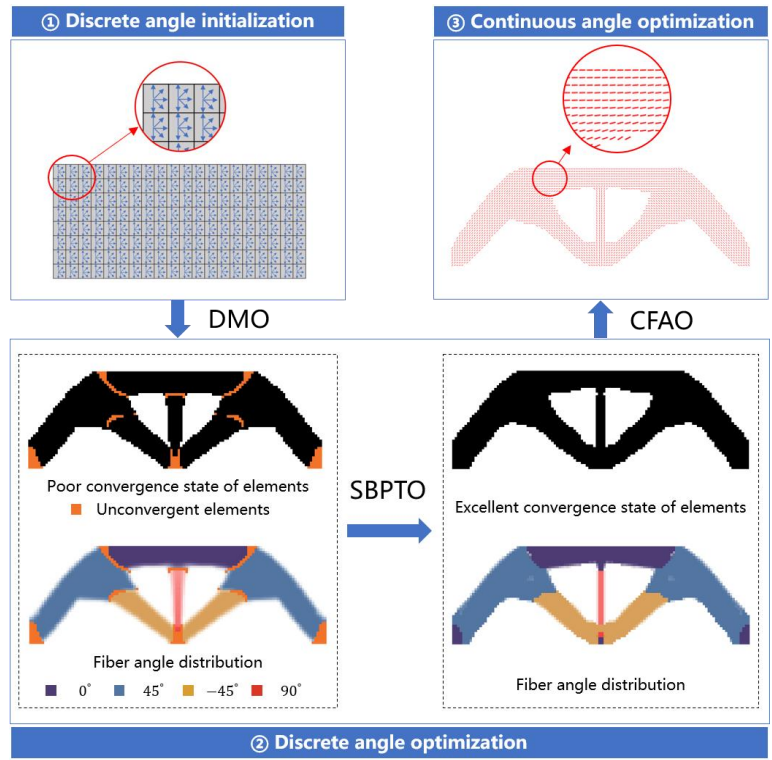

**Fig. 1**. Framework of the DSCO

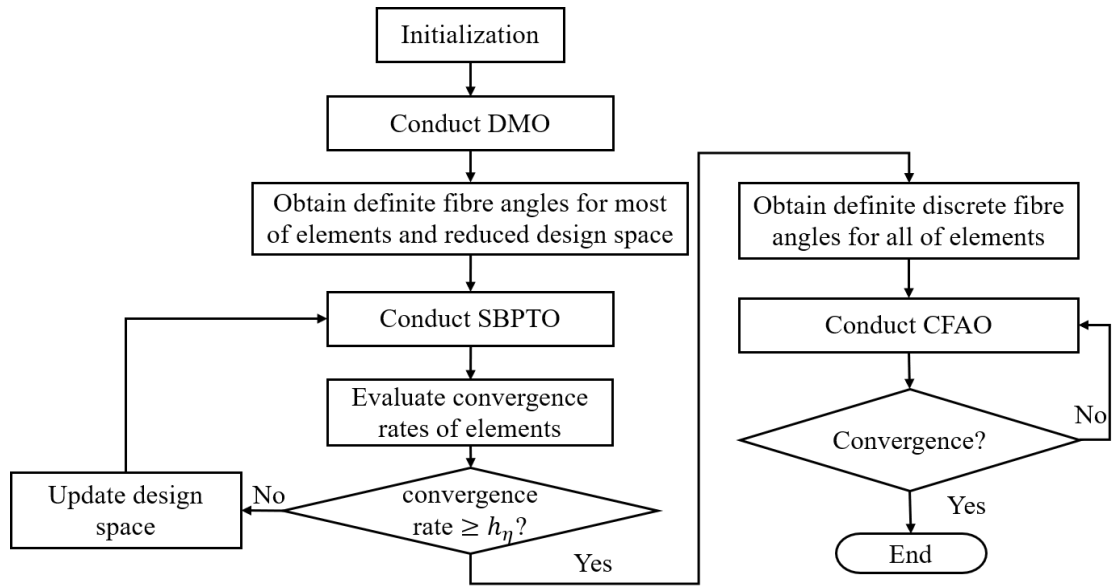

**Fig. 2.** The flow chart of the concurrent optimization

## 2.1. Topology optimization formulation

In this study, the objective is to minimize the compliance of structure by the concurrent optimization of topology and fibre angle orientation. It can be updated as follows:

$$\begin{aligned}
&find : (\boldsymbol{\chi}_1, \boldsymbol{\chi}_2, \ldots \boldsymbol{\chi}_n, \boldsymbol{\alpha}_1, \boldsymbol{\alpha}_2, \ldots \boldsymbol{\alpha}_{n+1}, \boldsymbol{\rho}, \boldsymbol{\theta}) \\
&\min : c = \mathbf{U}^T \mathbf{K} \mathbf{U} \\
&s.t. : \mathbf{K}\mathbf{U} = \mathbf{F} \\
&\qquad V \leq f \times N \\
&\qquad \alpha_{i,j} \in [0,1];\ \chi_{i,j} \in [0,1];\ \rho_i \in [0,1]; \\
&\qquad \theta_i \in \left[ -\frac{\pi}{2}, \frac{\pi}{2} \right]
\end{aligned} \qquad (1)$$

where $\chi_{i,j}$ and $\alpha_{i,j}$ $(i=1,\ldots n)$ denotes the density of the $j$th candidate fibre angle in the $i$th element for DMO and SBPTO respectively. $\rho_i$ is the density of the $i$th element in the step of CFAO. $\theta_i$ is fibre orientation in the $i$th element. The objective, $c$, stands for compliance. $\mathbf{K}$ denotes the stiffness matrix. $\mathbf{U}$ denotes the displacement vector. $\mathbf{F}$ denotes the load vector respectively. $N$ is the total number of elements. $V$ is the material volume and $f$ is the allowed maximum volume.

## 2.2. Topology optimization with Discrete Angle Optimization

Firstly, the Discrete Material Optimization (DMO) is introduced. The purpose of the DMO method is utilized to couple the two geometrical scales and carry out macroscopic topology optimization and microscopic material selection simultaneously. The main idea of DMO is to minimize the objective function by selecting a material among a predefined set of materials. In this study, a predefined set of materials represents a set of fibre angles.

In the optimization process, the design variables $\chi_i$, the density of corresponding material, vary between 0 and 1. Therefore, the discrete material problem can be transformed into the continuous variable optimization problem. The $i$th element constitutive matrix, $\mathbf{D}_i^e$, can be expressed as a weighted sum of candidate constitutive matrixes corresponding to the predefined fibre angles, $\mathbf{D}_j$:

$$\mathbf{D}_i^e = \sum_{j=1}^n w_{i,j} \mathbf{D}_j^e = w_{i,1}\mathbf{D}_1 + w_{i,2}\mathbf{D}_2 + \ldots + w_{i,n}\mathbf{D}_n, \quad 0 \leq w_{i,j} \leq 1 \qquad (2)$$

where $n$ denotes the number of candidate fibre angles. The weight functions, $w_{i,j}$, are between 0 and 1. And 0 implies giving up a material while 1 means choosing a material in physical sense. In this work, the weight function has been denoted as:

$$w_{i,j} = (\varepsilon + \chi_j^p)\Pi_{k=1}^{n}(\varepsilon + (1-\chi_{k\neq j}^p)) \tag{3}$$

where $\chi_i (i=1,\ldots,n)$ denotes the density of the $i$th material. In order to be 'fair', the initial design variables, $\chi_i$, should be set uniformly between 0 and 1. $\varepsilon$ is a small positive number (for example $1\times 10^{-9}$). The power, $p$, is a penalty for the intermediate values $\chi_i$ which is aiming at pushing the design variables to 0 or 1.

In this study, one candidate material corresponds to a candidate discrete fibre angle. From Eqs. (2-3), increasing the number of candidate materials tends to increase design variables. At the same time, increasing the number of elements or the number of laminates will significantly increase the number of design variables. The computation cost and storage of the optimization would also be increased rapidly.

In the DMO, the design variables can be updated by using MMA solver which is a gradient-based solver. Combined with Eq. (1), the derivative of compliance with respect to $\chi_{i,j}$ can be given by

$$\frac{\partial c_i}{\partial \chi_{i,j}} = -\mathbf{u}_e^T \frac{\partial \mathbf{k}_e}{\partial \chi_{i,j}} \mathbf{u}_e \tag{4}$$

where $c_i$ is the compliance of the $i$th element. $\mathbf{u}_e$ is the displacement component vector in the $i$th element. $\mathbf{k}_e$ is the element stiffness matrix.

$\dfrac{\partial \mathbf{k}_e}{\partial \chi_{i,j}}$ is obtained by

$$\frac{\partial \mathbf{k}_e}{\partial \chi_{i,j}} = \int_{\Omega_e} \mathbf{B}^T \frac{\partial \overline{\mathbf{D}}_i^e}{\partial \chi_{i,j}} \mathbf{B} d\Omega = \sum_{k=1}^{n} \int_{\Omega_e} \mathbf{B}^T \lambda_k \frac{\partial w_{i,k}}{\partial \chi_{i,j}} \mathbf{D}^e \lambda_k^T \mathbf{B} d\Omega \tag{5}$$

where

$$\mathbf{D}^e = \begin{bmatrix} D_{11}^e & D_{12}^e & 0 \\ D_{21}^e & D_{22}^e & 0 \\ 0 & 0 & D_{33}^e \end{bmatrix} \tag{6}$$

$$\boldsymbol{\lambda}_k = \begin{bmatrix} \cos^2\theta_k & \sin^2\theta_k & -\sin 2\theta_k \\ \sin^2\theta_k & \cos^2\theta_k & \sin 2\theta_k \\ \sin\theta_k \cos\theta_k & -\sin\theta_k \cos\theta_k & \cos^2\theta_k - \sin^2\theta_k \end{bmatrix} \tag{7}$$

Substituting Eqs. (2) and (3) into $\dfrac{\partial w_{i,k}}{\partial \chi_{i,j}}$,

$$\frac{\partial w_{i,k}}{\partial \chi_{i,j}} = \begin{cases} p\chi_j^{p-1} \prod_{k=1}^{n}\left(\varepsilon + \left(1 - \chi_{k\neq j}^p\right)\right), & \text{if } k = j \\ \left(-p\chi_j^{p-1}\right)\left(\varepsilon + \chi_m^p\right)\prod_{k=1}^{n}\left(\varepsilon + \left(1 - \chi_{k\neq m, k\neq j}^p\right)\right), & \text{otherwise} \end{cases} \tag{8}$$

For all numerical examples, the convergence criterion for the step of the DMO is as follow,

$$\varepsilon_i = \frac{\sum_{j=i-4}^{i}\left(c(j) - \bar{c}\right)^2}{5} \leq \varepsilon_0 \tag{9}$$

where $\varepsilon_i$ refers to the variance of the last five compliance values of the $i$ th iteration in DMO. $\varepsilon_0$ is set as a small constant.

**Algorithm 1:** Discrete material optimization (DMO)
**Step 1**: Set design domain, objective volume fraction, filter radius.
**Step 2**: each element is assigned a set of predefined candidate angles, $\chi_i (i=1,\ldots n)$, which are set uniformly.
**Step 3**: Calculate element constitutive matrix by Eq. (2) and carry out FEA.
**Step 4**: Calculate the compliance and the sensitivities using Eq. (5).
**Step 5**: Apply the MMA to determine the $n \times N_e$ volume fractions ($N_e$ is the number of elements).
**Step 6**: if the convergence criterion in Eq. (9) is satisfied, continue. If not, return to Step 3.

In order to elaborate the process of optimization better, it is necessary to define the term "fibre convergence". Ideally, in the process of discrete fibre angle optimization, the artificial density of one material is pushed to 1 as well as the densities of other materials are pushed to 0, or all the densities of materials are pushed to 0 for each element. However, there are some elements that cannot meet the above requirements. Several numerical examples in [39] indicates that convergence rate in DMO is limited and the excellent fibre convergence cannot be achieved. Therefore, it is required to use fibre convergence to evaluate the proportion of elements that meet the above criteria.

The element which satisfies the inequality Eq. (10) could be regarded as a converged one. In other words, the element selects a definite fibre angle from candidate angles.

$$\chi_i \geq \eta \sqrt{\chi_1^2 + \chi_2^2 + \ldots + \chi_p^2} \tag{10}$$

where $\eta$ indicates the tolerance level and it is typically in the interval [0.95, 0.995].

Fibre convergence, the ratio of converged elements in total elements, $h_\eta$, is depicted as:

$$h_\eta = \frac{N_c^e}{N^e} \tag{11}$$

where $h_{0.95} = 1$ simply means that for every element, there is one single volume fraction contributing more than 99.5% to the Euclidian norm of volume fractions.

Due to the elements located in "fuzzy regions", it is difficult to achieve excellent fibre convergence in DMO. These "fuzzy regions" are usually located at the junction of different angles or near where constraints and forces are applied. In these regions, the force transfer path is commonly complex. Generally, the sensitivities of the objective function to each fibre angle variable of the element located in these "fuzzy regions" are similar. This makes it difficult for those elements to make a choice in the candidate fibre angles in the stage of DMO. To handle the convergence problem, the sequential binary-phase topology optimization (SBPTO) with good selection performance in multi-phase is thus introduced. Multi-phase optimization is decomposed into a series of two-phase topology optimization.

Multi-phase with (n+1)-phase represents $n$ candidate fibre angles which have been defined in DMO and one void-material. The SBPTO decomposes the (n+1)-phase topology optimization method into a series of binary phases topology optimization.

After cyclic calculation, one candidate material would be selected while the remained materials have been discarded.

For the convenience of representation, we will represent the aforementioned design variables as $\alpha_i (i=1,\ldots,n+1)$ which determine the distribution of the material corresponding to the $i$th phase. In each element, the design variables should be equal to unity.

$$\sum_{i=1}^{n+1} \alpha_i = 1, \quad 0 \leq l_i \leq \alpha_i \leq u_i \leq 1 \tag{12}$$

where $l_i$ and $u_i$ represent the upper and lower boundaries respectively.

$n(n+1)$ binary phase topology sub-problems are involved in the SBPTO with $(n+1)$-phase. During the solution of each binary phase subproblem, $n-1$ phases in the unfrozen element should be fixed while the remained two phases would be active phases. The two active phases in unfrozen elements are presented by $a$ and $b$, and volume fractions are denoted by $\alpha_a^u$ and $\alpha_b^u$ respectively. And the fixed phases are denoted by $\alpha_i^u (i \neq \{a,b\})$. The sum of the two active phases $r_{ab}^u$ can be calculated as:

$$r_{ab}^u = 1 - \sum_{\substack{i=1 \\ i \neq \{a,b\}}}^{n+1} \alpha_i^u \tag{13}$$

The elements whose design variables satisfies Eq. (14) would be fixed in the next binary-phase topology optimization. The design variables of unfrozen elements form a reduced design space. This reduced design space serves as the initial design space for binary-phase topology optimization.

$$\alpha_a > \lambda \quad \text{or} \quad \alpha_b > \lambda \quad \text{or} \quad \begin{cases} \alpha_a = 0 \\ \alpha_b = 0 \end{cases} \tag{14}$$

where satisfying the first or second inequality means that the fibre angle of element is determined and satisfying the third inequality means both of a and b phase in the element are not selected. Therefore, the design variables in these elements that satisfy the inequality are fixed in this binary-phase topology optimization. $\lambda$ is set close to 1 but less than 1. In this work, we set $\lambda$ as 0.99.

For each binary-phase topology sub-problem, $\alpha_a^u$, is taken as the design variable.

The optimization of the binary-phase topology sub-problem refers to the SIMP method in [44]. After sub-problem optimization, the artificial densities of $b$, $\alpha_b^u$, would be computed as follow:

$$\alpha_b^u = r_{ab}^u - \alpha_a^u \tag{15}$$

the upper bound of phase $a$, $u_{a,temp}^u$, in each iteration should satisfy:

$$u_{a,temp}^u = \min\left(u_a^u, r_{ab}^u\right) \tag{16}$$

In the sequential binary-phase topology optimization, the derivative of compliance with respect to $\alpha_{i,a}$ can be obtained:

$$\frac{\partial c_i}{\partial \alpha_{i,a}} = -\mathbf{u}_e^T \frac{\partial \mathbf{k}_e}{\partial \alpha_{i,a}} \mathbf{u}_e \tag{17}$$

where

$$\frac{\partial \mathbf{k}_e}{\partial \alpha_{i,a}} = \int_{\Omega_e} \mathbf{B}^T \frac{\partial \overline{D_e}}{\partial \alpha_{i,a}} \mathbf{B} d\Omega = \int_{\Omega_e} \mathbf{B}^T \lambda \frac{\partial D_e}{\partial \alpha_{i,a}} \lambda^T \mathbf{B} d\Omega \tag{18}$$

The alternating active-phase algorithm can be depicted as:

---
**Algorithm 2:** Sequential Binary-Phase Topology Optimization

$$\boldsymbol{\alpha}^u = \left(\alpha_1^u, \ldots, \alpha_{n+1}^u\right) = \left(\chi_1, \chi_2, \ldots, \chi_n, 1 - \sum_{i=1}^n \chi_i\right)$$

**Repeat**
    **for** a=1 **to** p **do**
        **for** b=1 **to** p (except for a) **do**
            $\boldsymbol{\alpha}^u$ = solution of binary phase subproblem
        **end**
    **end**
    update filter_params
    iter_out = iter_out + 1
**Until** (iter_out>threshold) **or** ($h_\eta$> criteria)

---

The flow chart of the SBPTO is shown as **Fig. 3**:

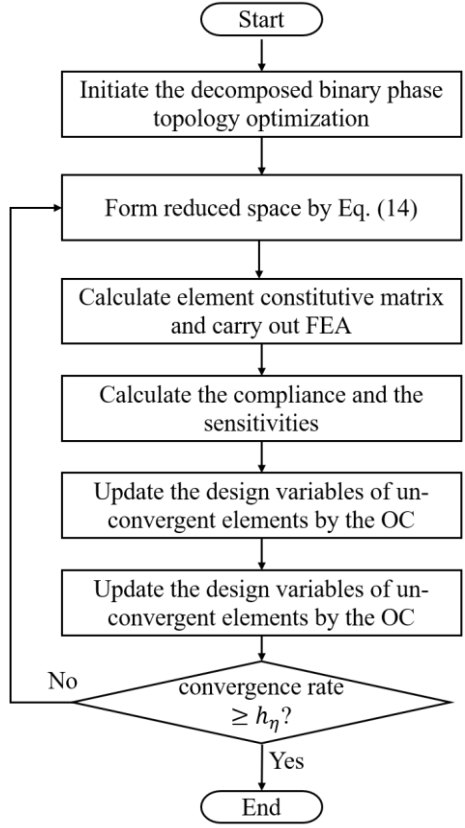

**Fig. 3**. The flow chart of the SBPTO

## 2.3. Representation of fibre angle

To improve the fibre continuity, a continuous interpolation function is introduced. It is a fibre continuity filter which is a spatial filter [42]. In the design domain, each fibre angle of design point $\theta_i$, $i=1,\ldots n$ at $p_i=(x_i,y_i,z_i)$ can be interpolated by these design variables in a "influence area", that satisfies the condition, $\|\mathbf{x}-\mathbf{p}_i\| \leq R_c$, where $R_c$ is called the cut-off radius. The interpolation function is defined as:

$$\Theta(\mathbf{x}) = \sum_{i \in I_x} w_i(\mathbf{x})\theta_i, \quad \|\mathbf{x}-\mathbf{p}_i\| \leq R_c \tag{19}$$

the fibre angle, $\Theta(\mathbf{x})$, should be constrained as $\Theta(\mathbf{x}) \in [-90°, 90°]$. The initial fibre angle configuration is determined by the results of discrete optimization above. And $w_i(\mathbf{x})$ is the weighting function which is defined as follow.

$$w_i(\mathbf{x}) = \frac{H_{ei}\alpha_i}{\sum_{j \in I_x} H_{ej}}$$
$$H_{ei} = \max(0, r_{min} - \|\mathbf{x}-\mathbf{p}_i\|) \tag{20}$$

where $\alpha_i$ is the volume fraction.

In the continuous fibre angle optimization, the derivative of compliance with respect to the design variable $\rho_i$ can be given by

$$\frac{\partial c_i}{\partial \rho_i} = -\mathbf{u}_e^T \frac{\partial \mathbf{k}_e}{\partial \rho_i} \mathbf{u}_e \tag{21}$$

Substituting Eq. (16) into $\frac{\partial \mathbf{k}_e}{\partial \rho_i}$

$$\frac{\partial \mathbf{k}_e}{\partial \rho_i} = \int_{\Omega_e} \mathbf{B}^T \frac{\partial \overline{D}_e}{\partial \rho_i} \mathbf{B} d\Omega = \int_{\Omega_e} \mathbf{B}^T \lambda \frac{\partial D_e}{\partial \rho_i} \lambda^T \mathbf{B} d\Omega \tag{22}$$

The derivative of compliance with respect to $\theta_i$ can be developed

$$\frac{\partial \mathbf{k}_e}{\partial \theta_i} = \int_{\Omega_e} \mathbf{B}^T \frac{\partial \overline{D}_e}{\partial \theta_i} \mathbf{B} d\Omega = \int_{\Omega_e} \mathbf{B}^T \frac{\partial \lambda}{\partial \theta_i} D_e \lambda^T \mathbf{B} d\Omega + \int_{\Omega_e} \mathbf{B}^T \lambda D_e \frac{\partial \lambda^T}{\partial \theta_i} \mathbf{B} d\Omega \tag{23}$$

where

$$\frac{\partial \lambda}{\partial \theta_i} = \begin{bmatrix} -2\sin\theta_i \cos\theta & 2\sin\theta_i \cos\theta_i & -2\cos 2\theta_i \\ 2\sin\theta_i \cos\theta_i & -2\sin\theta_i \cos\theta_i & 2\cos 2\theta_i \\ \cos 2\theta_i & -2\cos 2\theta_i & -4\sin\theta_i \cos\theta_i \end{bmatrix} \tag{24}$$

## 3. Numerical examples

This section demonstrates the validity of the proposed method with four numerical examples: an MBB beam with three in-plane loads, an L-shape beam with one in-plane load, a cantilever beam with one in-plane load and a cantilever beam with multiple loads. For simplicity, all of examples employ uniform meshes which size are $1 \times 1$. And in the step of the SBPTO, the convergence criterion, $h_\eta$, is set as 0.99.

### 3.1. MBB beam

The first numerical example is a 2D MBB beam structure. The boundary condition is shown in Fig. 4. The points at the one-fourth and three-fourth of the top edge are applied external force $F = 1$. The middle point of bottom edge is applied external force $2F$. The rectangular domain size is $120 \times 40$ and is meshed by 4-nodes regular quadrilateral elements. The material parameters with orthotropic properties are given as $D_{11}^e = 0.5448$, $D_{12}^e = 0.0383$, $D_{22}^e = 0.1277$, $D_{33}^e = 0.0456$. Only one layer has been used. The desired volume fraction is 0.5. the minimum filter radius $r_{min}$ is 1.5. For comparison, these cases are investigated by the CFAO and DSCO separately. In the cases(a)~(d) of the CFAO, the initial fibre angles of CFAO are set as case(a)~(d) and candidate angles sets of DSCO are set as cases(e)~(n) in Table 1.

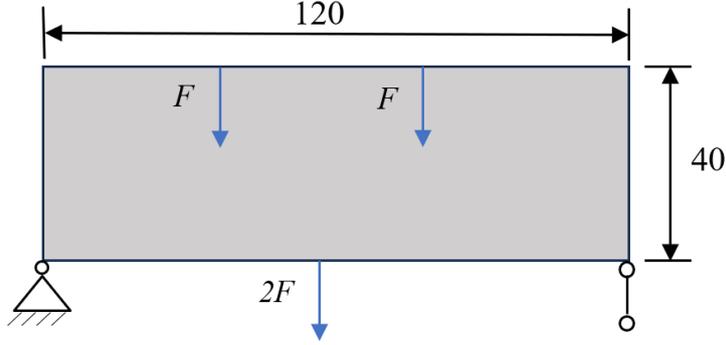

**Fig. 4.** The MBB beam.

**Table 1.** the angle setting of the CFAO and DSCO

|  | case | $\theta_{initial}^e$ | case | $\theta_{initial}^e$ |
|---|---|---|---|---|
| CFAO | case(a) | [0°] | case(b) | [90°] |
|  | case(c) | [45°] | case(d) | [−45°] |
| DSCO | case(e) | [0°, 90°] | case(f) | [0°, −30°, 30°, 90°] |
|  | case(g) | [0°, −60°, 60°, 90°] | case(h) | [0°, −45°, 45°, 90°] |
|  | case(i) | [0°, −45°, 45°, 90°, −30°, 30°] | case(j) | [0°, −45°, 45°, 90°, −60°, 60°] |
|  | case(k) | [0°, −45°, 45°, 90°, −30°, 60°] | case(l) | [0°, −45°, 45°, 90°, 30°, 60°] |
|  | case(m) | [0°, −45°, 45°, 90°, 30°, −60°] | case(n) | [0°, −45°, 45°, 90°, −30°, −60°] |

The optimized results are shown in Fig. 5 and Fig. 6. Fig. 5 shows that the iteration number of CFAO is generally lower than that of DSCO, but the compliance of CFAO is particularly sensitive to the initial fibre angle setting. The maximum and minimum compliance in cases (a)~(d) are 622.42 and 328.32 respectively. That is consistent with the previous introduction in section 1. The optimized compliances of the DSCO are shown in Fig. 5(e)~(n). It can be seen that compared with CFAO, the optimized optimization results of the DSCO are more stable for different initial fibre configurations especially for the cases(f)~(n). Since there are only 2 angles in the initial candidate angle set, the compliance of the case(e) has the largest fluctuation in DSCO. That indicates that the optimization ability of the candidate angle set containing two angles is poor. The initial angle sets of the rest cases of DSCO are 4 angles or 6 angles. It can be found that the compliance does not decrease with the increase of the number of angles in the initial angle set and the case(h) with initial angle set $[0°,-45°,45°,90°]$ has minimum compliance 302.7. At the same time, we can also see that in cases(e)~(n), the number of iterations increases significantly as the number of angles of the angle set increases. For the case(h), the compliance is 302.7 and the iteration number is 434. Meanwhile for the case(d), which has the minimum compliance in the cases of CFAO, the compliance is 328.32 and the iteration number is 350. Compared with case(d), the compliance of case(h) decreases by 7.8%, which is crucial for the optimization, even with higher computational costs. It is also worth mentioning that since CFAO is very sensitive to initial fibre angle setting, we cannot select an appropriate initial fibre angle setting directly. To obtain good result, we need to conduct several CFAO with different initial fibre angle settings, which would lead to the increase of computing costs substantially.

Although we cannot prove that the optimized solution of case(h) is a globally optimal solution, we can conclude that in this MBB beam example, the DSCO does reduce the risk of getting stuck in a local optimum and the initial angle set $[0,-45°,45°,90°]$ in DSCO provides good result.

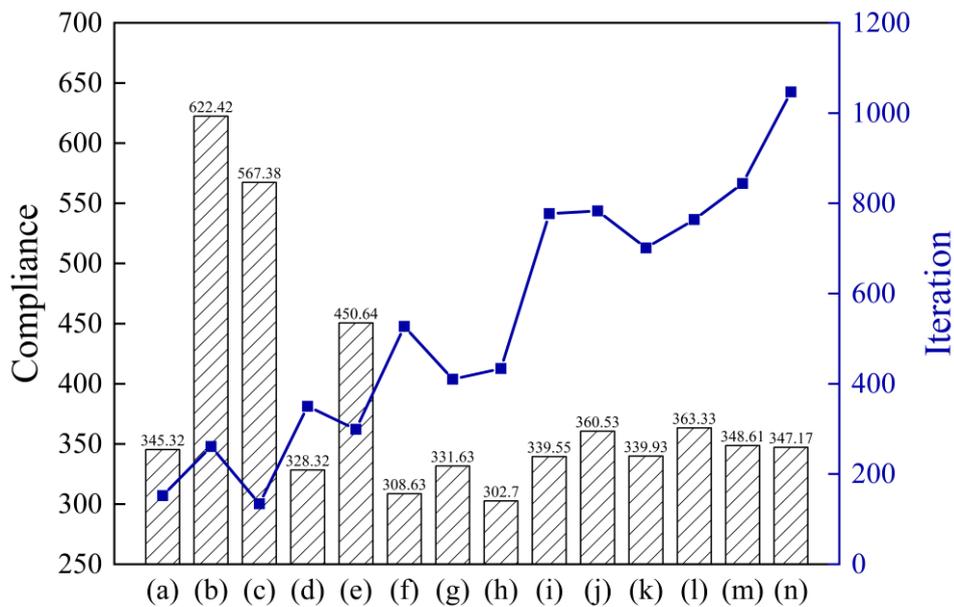

**Fig. 5**. The compliance and iteration of case(a)~(l)

In the optimized fibre layout figures of Fig. 6, the short red lines represent the fibre orientation of each element. As shown in Fig. 6, the fibre continuity in most of design domain is good. The few elements with poor continuity are mainly concentrated around the nodes where loads or constraints are applied or some elements in the area of rapid geometric change. For the overall structure, different initial angle settings always result in distinctively different structures. Due to the symmetry of the structure, half of the optimized fibre layout of case(h) is showed, see Fig. 7. It can be seen that fibre orientation of few elements located near the load varies erratically.

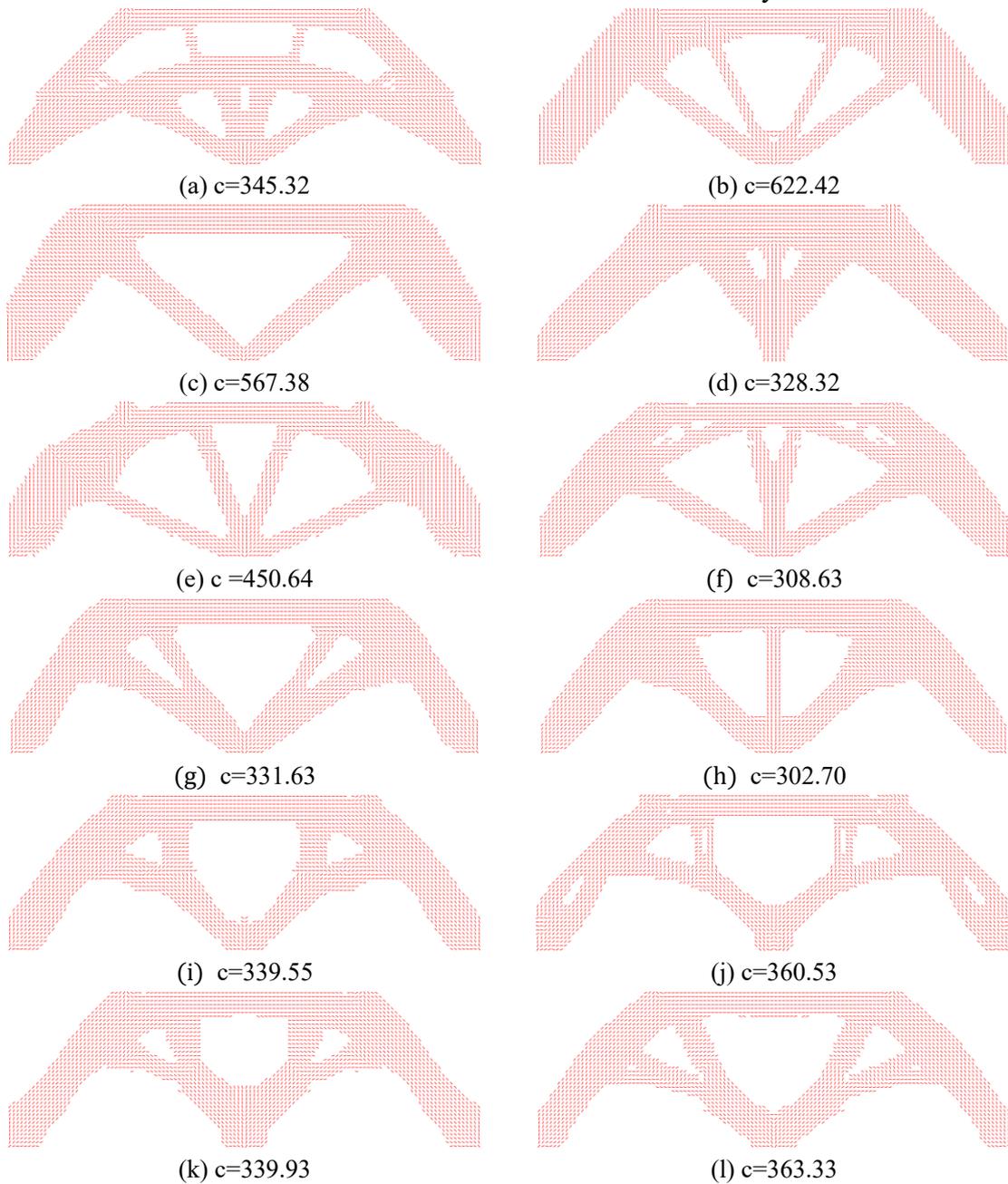

(a) c=345.32  (b) c=622.42
(c) c=567.38  (d) c=328.32
(e) c =450.64  (f) c=308.63
(g) c=331.63  (h) c=302.70
(i) c=339.55  (j) c=360.53
(k) c=339.93  (l) c=363.33

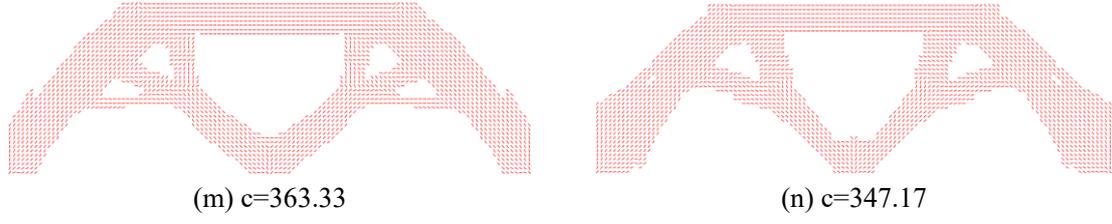

(m) c=363.33  (n) c=347.17

**Fig. 6.** The optimized fibre layout of CFAO and DSCO for MBB beam.

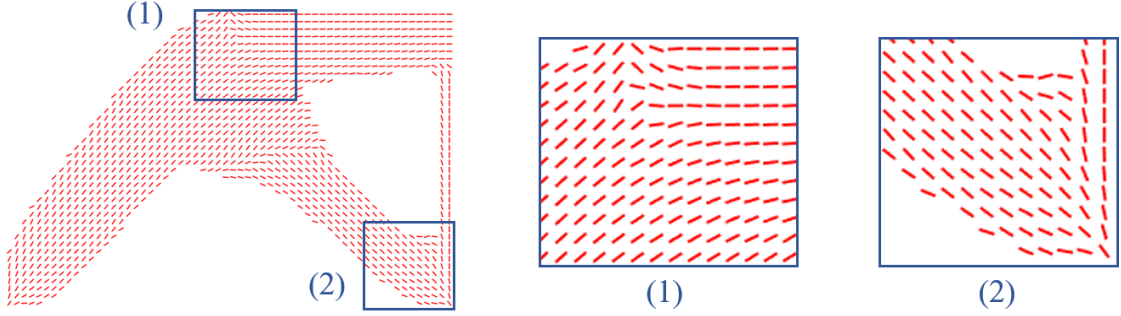

**Fig. 7**. Half of the optimized fibre layout of case(h)

Fig. 8 shows the objective function and convergence rate histories of case(h) with initial fibre angle set $\theta^e_{initial} = [0°, -45°, 45°, 90°]$ by using the DSCO for MBB beam. In the optimized structure figures, orange areas represent unconvergent elements and areas with other different colours represent different fibre angles. The total iteration of case(h) is 434 and the optimized compliance reaches a low value 302.70. Since the weighting function is not normalized, the compliance would be unrealistically high initially and this phenomenon has no effect on the final result [37]. The first step, DMO, stops when the convergence condition is met at the 60th iteration. The convergence rate, $h_{0.95}$, reaches 0.96 at the end of DMO. The second step, SBPTO, lasts 305 iterations. The convergence rate, $h_{0.95}$, reaches 1 at the end of SBPTO. The first two steps determine the initial values $\theta_i$ of the CFAO which realizes the continuous fibre angle optimization. It is worth mentioning that the orientation of the fibre is mostly parallel to the member direction.

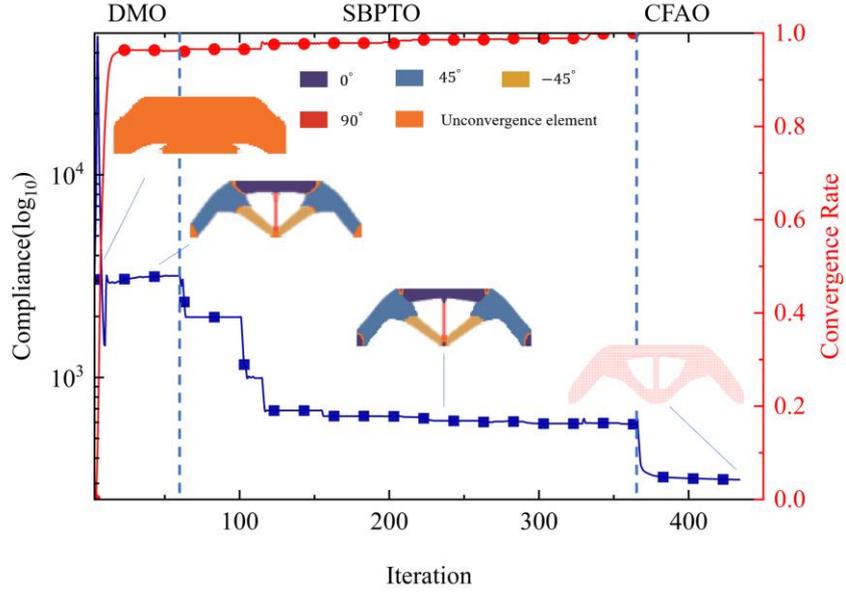

**Fig. 8**. Histories of objective and convergence rate with the DSCO for MBB beam

## 3.2. L-shape beam

The second numerical example is a 2D L-shape beam. The boundary condition and the size of the 2D L-shape structure are shown in Fig. 9. The top point of the right edge is applied an external force $F=1$. The convergence criterion, $\varepsilon_0$, is set as $10^{-2}$. The height and width of the beam are both 100. The beam is meshed by 4-nodes regular quadrilateral elements. The material parameters with orthotropic properties are given same with the case of MBB beam. Only one layer has been used. The desired volume fraction is 0.6 and the minimum filter radius $r_{min}$ is 1.5. The initial fibre angles of CFAO are and candidate angles sets of DSCO are shown in Table 2.

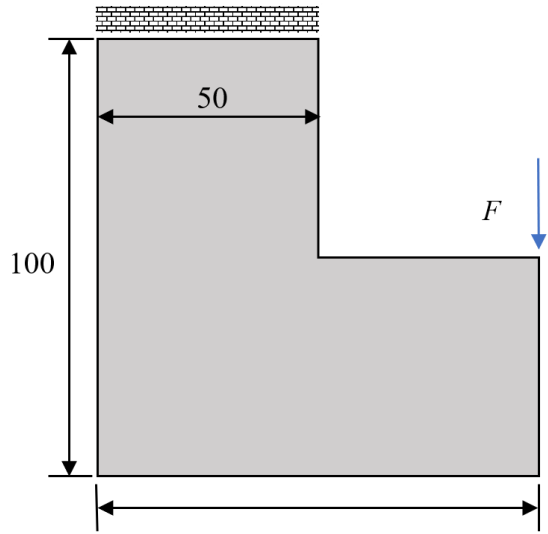

**Fig. 9.** L-shape beam

Table 2. The allocation of the initial fibre angle

| | case | $\theta^e_{initial}$ | case | $\theta^e_{initial}$ |
|---|---|---|---|---|
| CFAO | case(a) | [0°] | case(b) | [90°] |
| | case(c) | [45°] | case(d) | [−45°] |
| DSCO | case(e) | [0°, 90°] | case(f) | [0°, −30°, 30°, 90°] |
| | case(g) | [0°, −60°, 60°, 90°] | case(h) | [0°, −45°, 45°, 90°] |
| | case(i) | [0°, −45°, 45°, 90°, −30°, 30°] | case(j) | [0°, −45°, 45°, 90°, −60°, 60°] |

The optimized results are shown in Fig. 10 and Fig. 11. We can find that the iteration number of CFAO is generally lower than that of DSCO in Fig. 10. The minimum compliance of CFAO and DSAO are 188 and 170.58 respectively. Compared with CFAO, the optimized optimization results of the DSCO are more stable. The compliance of CFAO is particularly sensitive to the initial fibre angle setting. The maximum and minimum compliance in cases (a)~(d) are 326 and 188 respectively. The case(h) with initial angle set [0°, −45°, 45°, 90°] has minimum compliance 170.58. Meanwhile, we can also find that in cases(e)~(n), the number of iterations increases significantly as the number of angles of the angle set increases. For the case(h), the compliance is 170.58 and the iteration number is 354. Meanwhile for the case(d), which has the minimum compliance in the cases of CFAO, the compliance is 188 and the iteration number is 110. Compared with case(d), the compliance of case(h) decreases by 9.3%. Although it seems that the computation cost of CFAO is much lower, in order to get good result, the initial fibre angle settings are chosen based on trial and error method which will significantly increase the computation cost.

In the L-shape beam example, results show that the DSCO does reduce the risk of getting stuck in a local optimum and the initial angle set [0°, −45°, 45°, 90°] in DSCO provides good result.

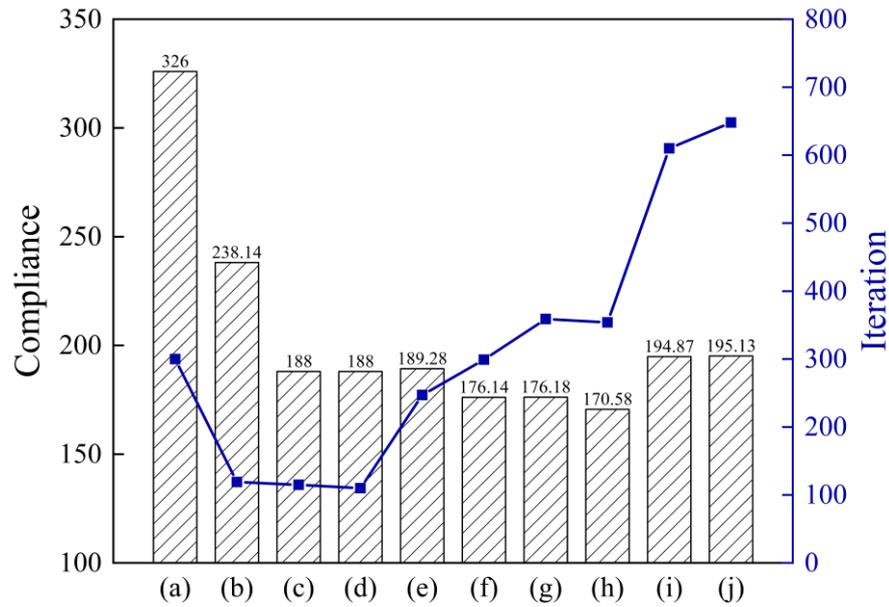

Fig. 10. The compliance and iteration of case(a)~(j)

As shown in Fig. 11, materials distribution of CFAO is relatively concentrated. The fibre continuity in most of design domain is also good. Few elements with poor continuity are mainly concentrated around the nodes where loads or constraints are applied and in the area of the corner of 'L' or low stress values. In order to show the fibre layout more clearly, optimized fibre layout of case(h) is showed in Fig. 12.

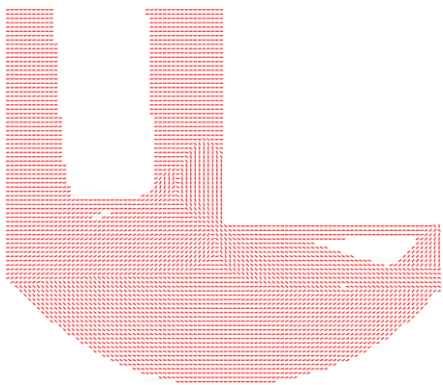     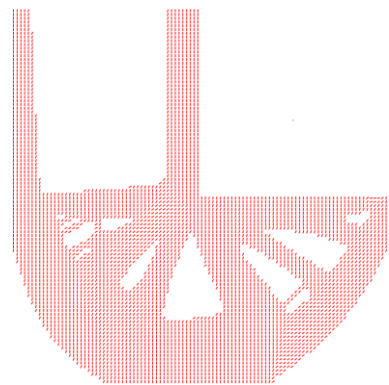

(a) c=326.66              (b) c=238.14

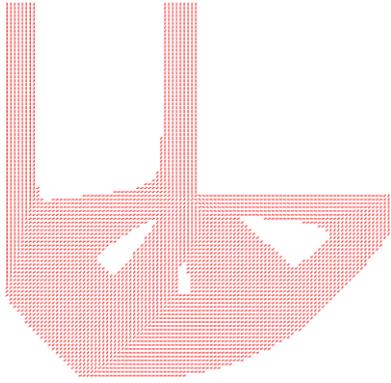

(c) c=188.10

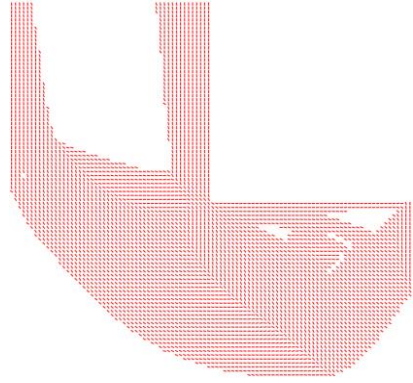

(d) c=188

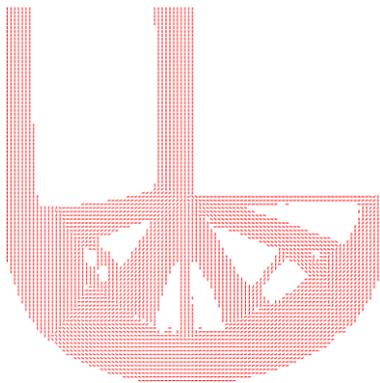

(e) c=189.28

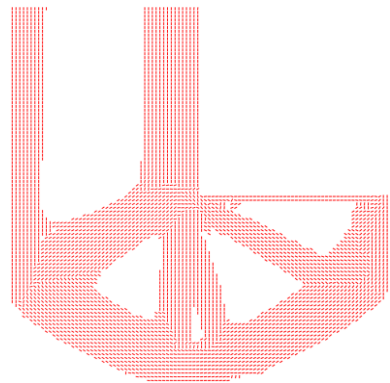

(f) c=176.14

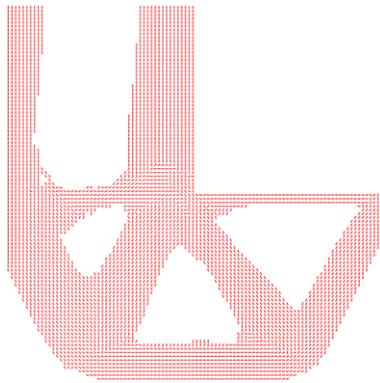

(g) c=176.18

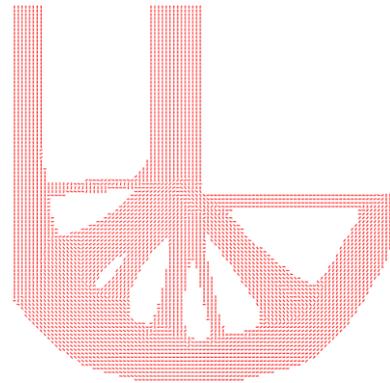

(h) c=170.58

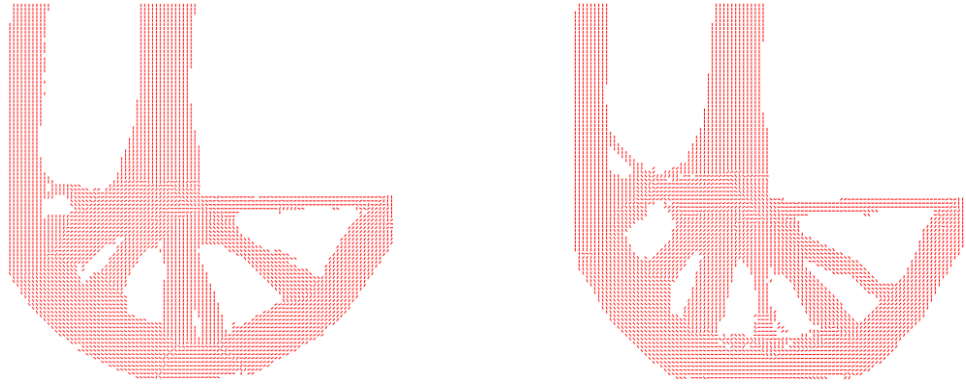

(i) c=194.87  (j) c=195.13

**Fig. 11**. The optimized results of CFAO and DSCO
for L-shape beam including fibre path

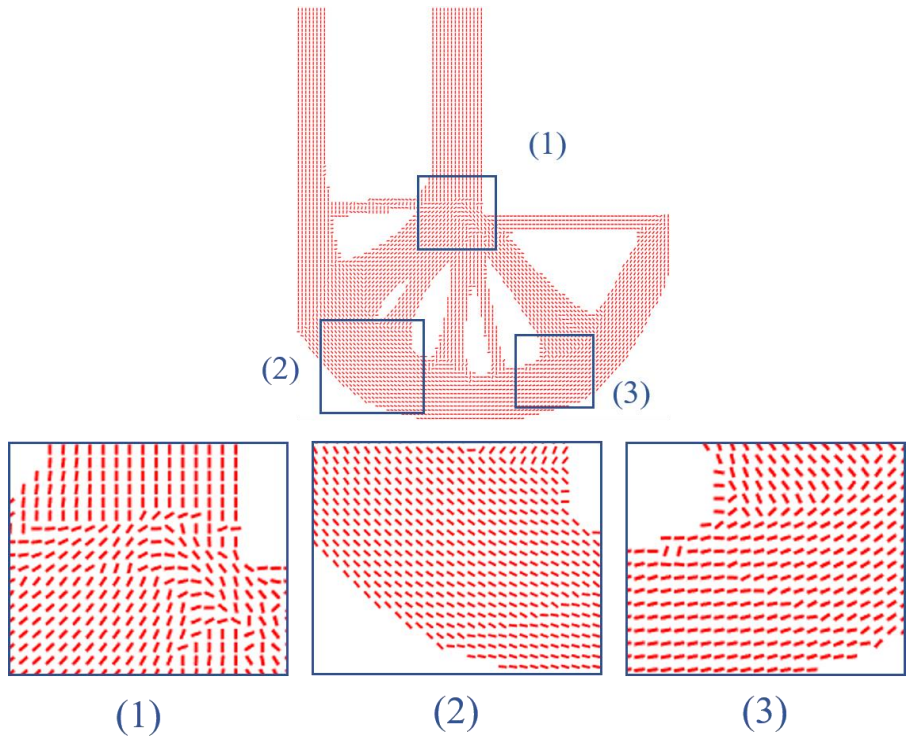

(1)  (2)  (3)

**Fig. 12.** the optimized fibre layout of case(h)

Fig. 13 shows the objective function and convergence rate histories of case(h) with initial fibre angle set $\theta_{initial}^e = [0°, -45°, 45°, 90°]$ by using the DSCO for MBB beam. The total iteration of case(h) is 354 and the optimized compliance reaches a low value 170.58. The first step, DMO, stops when the convergence condition is met at the 50th iteration. The convergence rate, $h_{0.95}$, reaches 0.95 at the end of DMO. The second step, SBPTO, lasts 246 iterations. The convergence rate, $h_{0.95}$, reaches 0.999 at the end of SBPTO. The third step, CFAO, lasts 58 iterations. We can find that the unconvergent elements in DMO mainly exists at the corner of the 'L' shape.

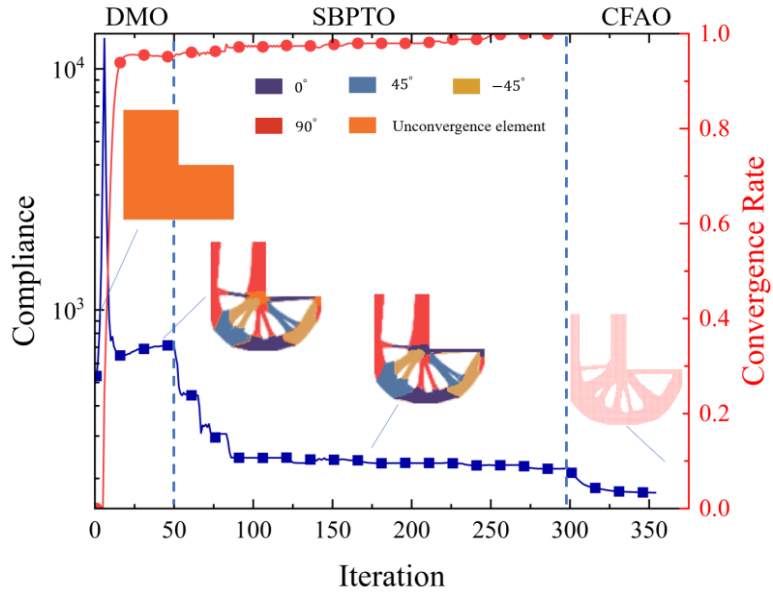

**Fig. 13.** Histories of objective functions and convergence rate with the DSCO for L-beam

### 3.3. Cantilever beam

**3.3.1 single load case**

The third numerical example is a 2D cantilever structure. The boundary condition is shown in Fig. 14. The left edge is fixed and the middle point of the right edge is applied an external force $F=1$. The convergence criterion, $\varepsilon_0$, is set as $10^{-2}$. The rectangular domain size is $50\times40$ and is meshed by 4-nodes regular quadrilateral elements. The material parameters with orthotropic properties are given as $E_x=2$, $E_y=1$, $G_{xy}=0.25$ and $v_{xy}=0.3$. Only one layer has been used. The desired volume fraction is 0.5. the minimum filter radius $r_{min}$ is 1.5. The allocation of the initial fibre angle sets of CFAO and DSCO are set as Table 2.

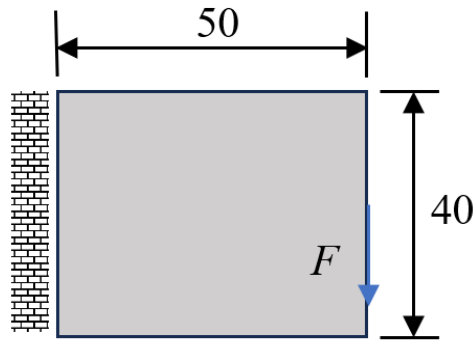

Fig. 14. The cantilever beam with single load.

The compliance and iteration number of case(a)~(j) are shown in Fig. 15. The cases(a)~(d) belong to CFAO and cases(e)~(j) belong to DSCO. DSCO's compliances are much more stable to the initial fibre angle settings than DSCO's. The compliance of case(h) with fibre angle setting $\theta_{initial}^{e} = [0°, -45°, 45°, 90°]$ is the minimum and the minimum compliance is 13.24. At the same time, the compliance of case(f) (13.25) is very close to the compliance of case(h) (13.24). In terms of the number of iterations, case(f) is 189 while case(h) is 224. Thus both initial fibre settings of case(f) and case(h) are good choices. Although we cannot prove that the optimized solution of case(h) is a globally optimal solution, we can conclude that the proposed method, DSCO, does reduce the risk of getting stuck in a local optimum.

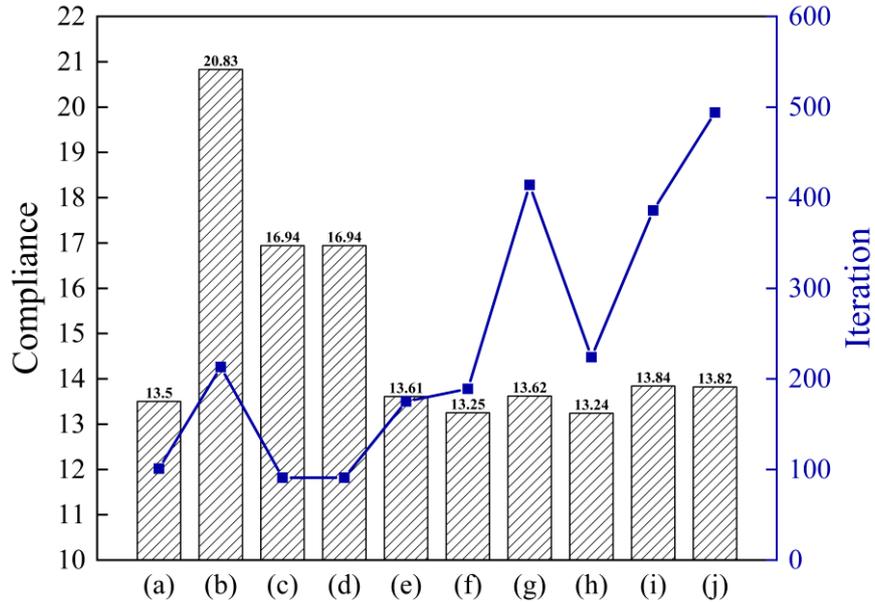

Fig. 15. The compliance and the iteration number of case(a)~(j).

The fibre optimized figures in **Fig. 16** show that the characteristics of the fibre continuity and the optimized structures with different initial angle settings are similar to those mentioned above. It is worth mentioning that the orientation of the fibre is mostly parallel to the member direction.

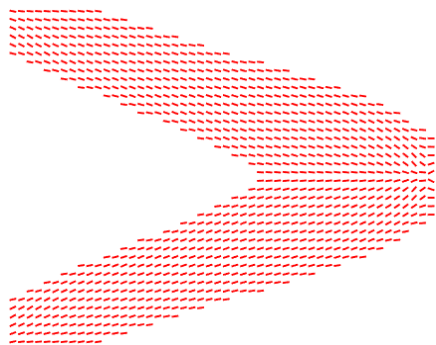

(a) c=13.50

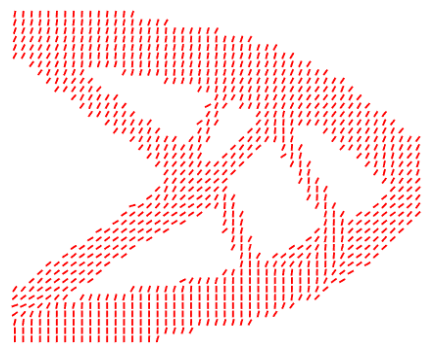

(b) c=20.83

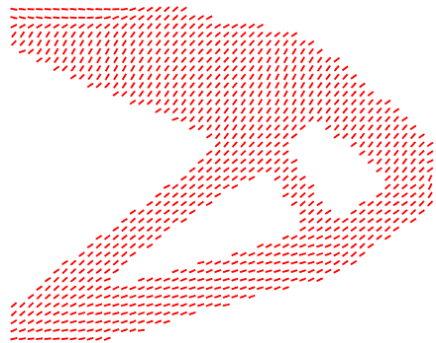

(c) c=16.94

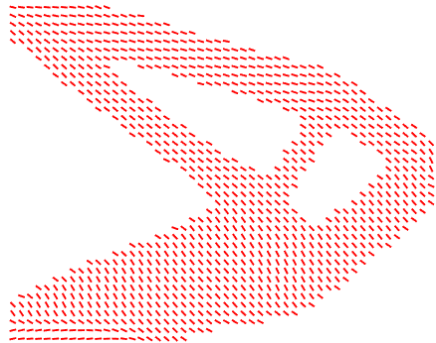

(d) c=16.94

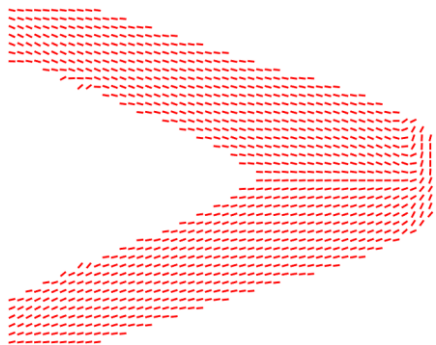

(e) c=13.61

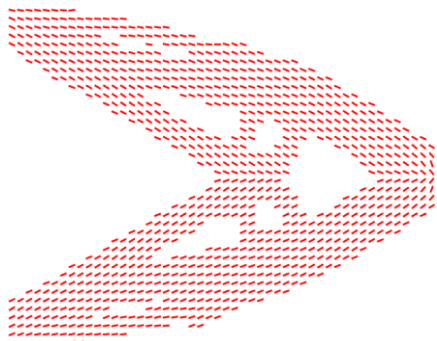

(f) c= 13.25

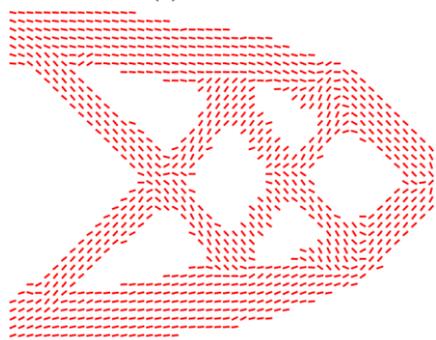

(g) c= 13.62

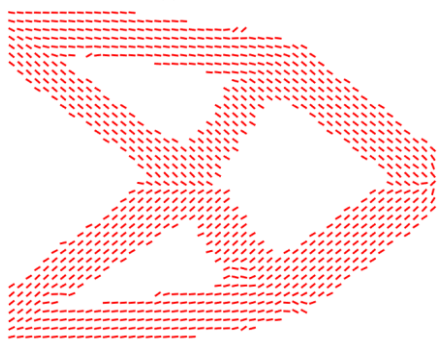

(h) c= 13.24

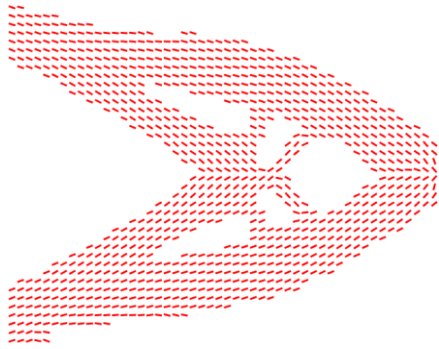
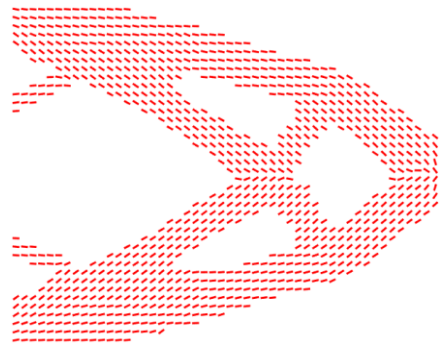

(i) c=13.84

(j) c=13.82

**Fig. 16.** The optimized fibre layout of CFAO and DSCO for Cantilever beam cases(a)~(j).

Fig. 17 shows the iterative process of the compliance and convergence rate of case(e)~(j). It can be found that the iterative processes are all stable. The convergence of each case can reach a very high level.

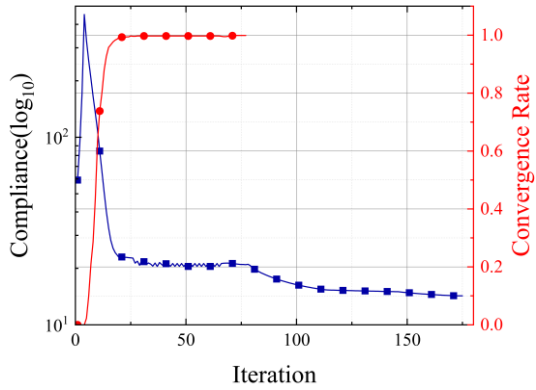
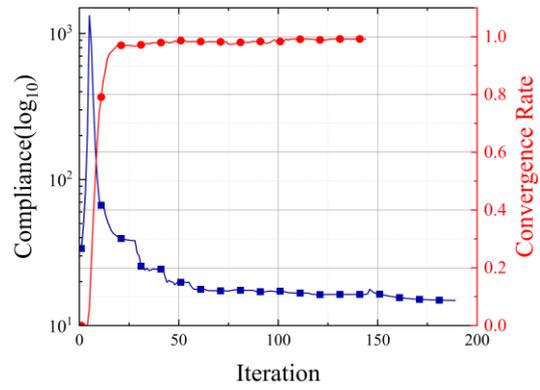

(e)　　　　　　　　　　　　　(f)

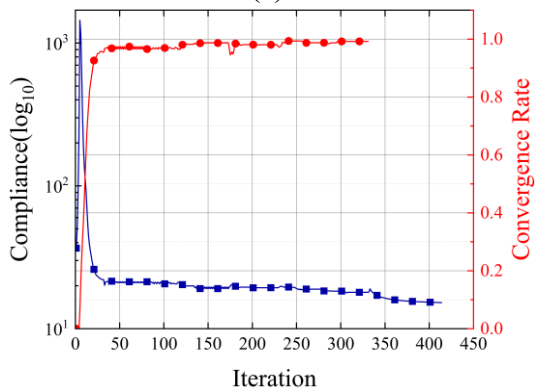
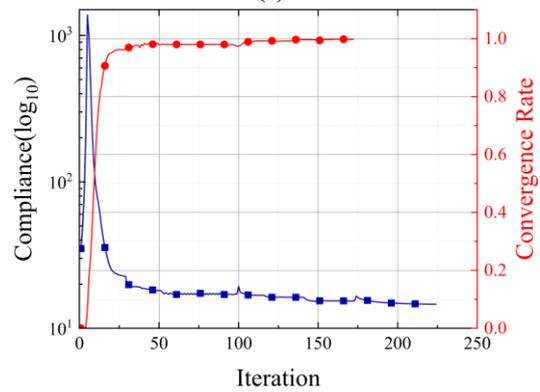

(g)　　　　　　　　　　　　　(h)
28

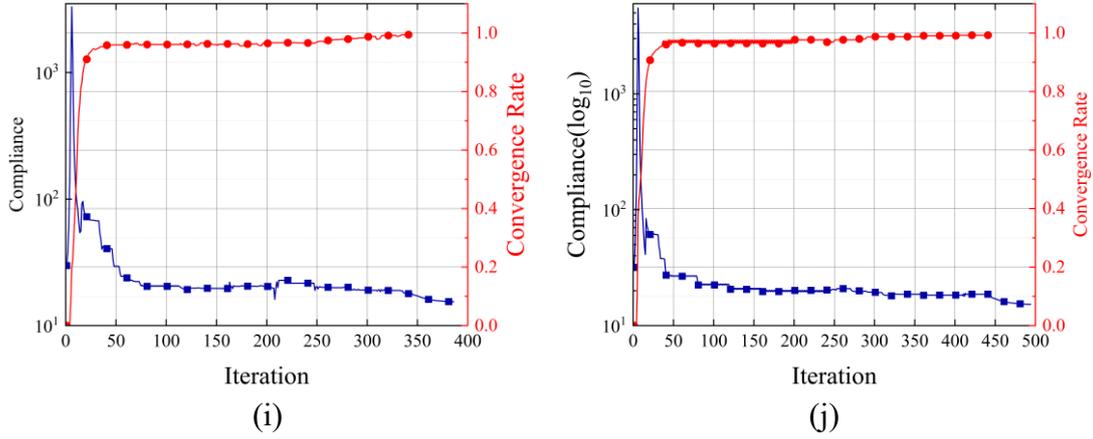

(i)                                        (j)

**Fig. 17**. The iterative histories for Cantilever beam cases(e)~(j) of DSCO

Fig. 18 shows the objective function and convergence rate histories of case(f) with initial fibre configuration $\theta_{initial}^{e} = [0°, -45°, 45°, 90°]$ by using the DSCO for Cantilever beam. The total iteration is 224 and the optimized compliance reaches a low value 13.24. The first step, DMO, stops when the convergence condition is met at the 29th iteration. The convergence rate, $h_{0.95}$, reaches 0.962 at the end of DMO. The second step, SBPTO, lasts 143 iterations. The convergence rate, $h_{0.95}$, reaches 0.9975 at the end of SBPTO.

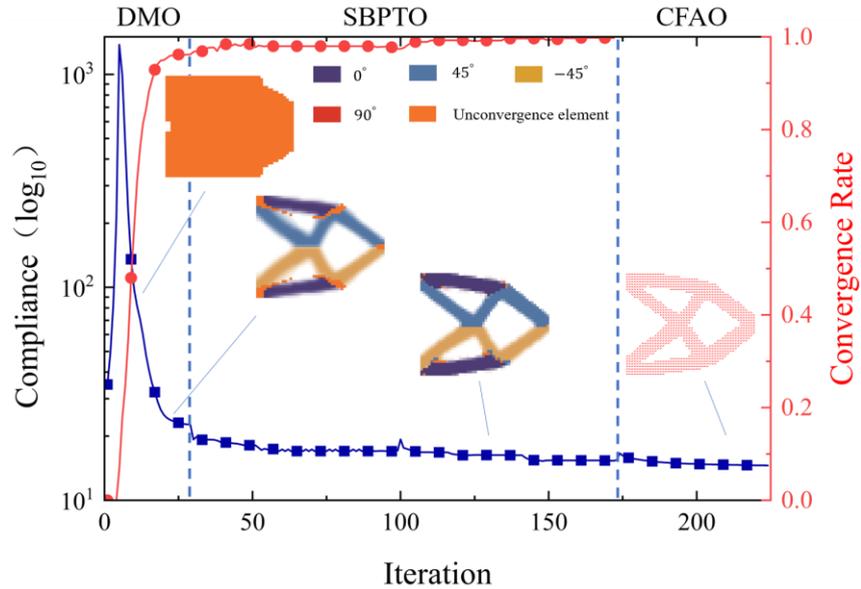

**Fig. 18**. Histories of objective function and convergence rate
of case(f) with the DSCO for Cantilever beam with single load.

### 3.3.2. muti-load case

The fourth numerical example is a 2D cantilever structure with multi-load. The boundary condition is shown in Fig. 19. The external force $F = 1$. The convergence criterion, $\varepsilon_0$, is set as $10^{-2}$. The rectangular domain size is $60\times40$ and is meshed by 4-nodes regular quadrilateral elements. The other parameters are same with the case in section 3.3.1.

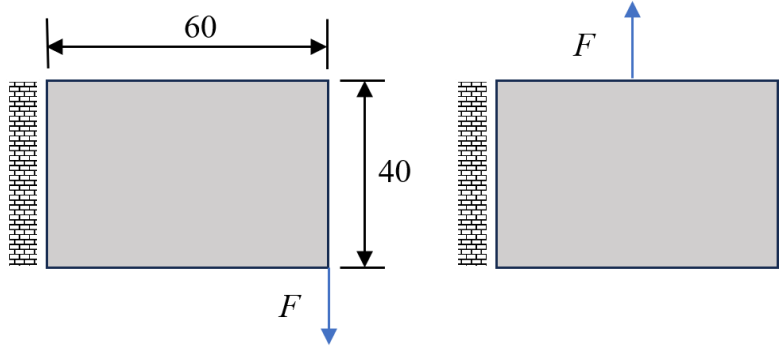

**Fig. 19**. The cantilever beam with multi-load.

The compliance and the iteration number of case(a)~(h) are shown in Fig. 20. The cases(a)~(d) are solved by CFAO while cases(e)~(h) use DSCO. The iteration number of CFAO is generally lower than that of DSCO, but the compliance of CFAO is particularly sensitive to the initial fibre angle setting. Except case(d), the compliances of CFAO are much higher than that of DSCO. In contrast, compliances of DSCO are much more stable for the initial fibre angle settings. The compliance of case(f) with fibre angle setting $\theta^e_{initial} = [0°, -45°, 45°, 90°]$ is the minimum (36.91). In the cases(e)~(f) of the DSCO, with the increase of the number of candidate angles, the iteration increases significantly. Although we cannot prove that the optimized solution of case(f) is a globally optimal solution, we can conclude that the proposed method, the DSCO, does reduce the risk of getting stuck in a local optimum.

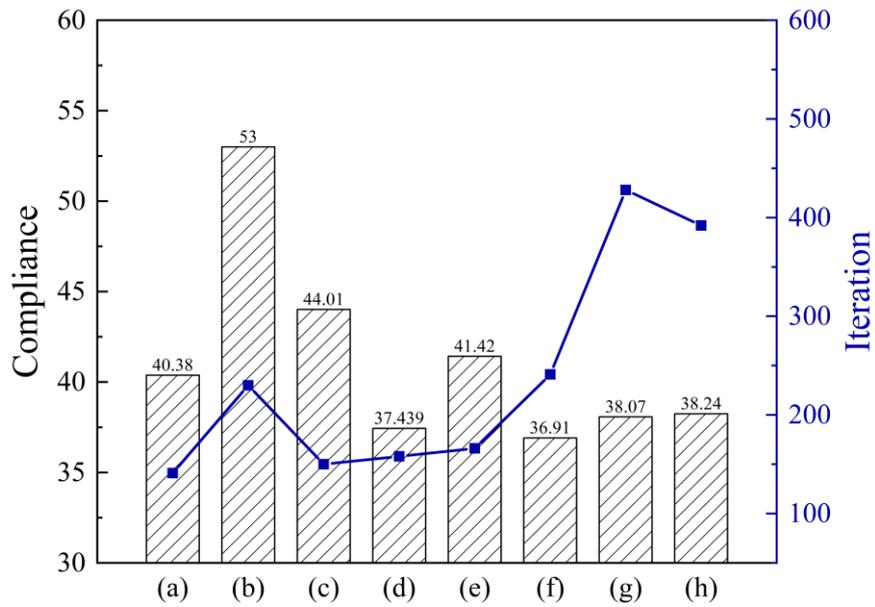
**Fig. 20.** The compliance and the iteration number of case(e)~(f).

In the fibre optimized figures of Fig. 21, the fibre continuity in most of design domain is good. The few elements with poor continuity are mainly concentrated around the nodes where loads or constraints are applied or some elements where there is a large variation in size. For the overall structure, different initial angle settings always result in very different structures.

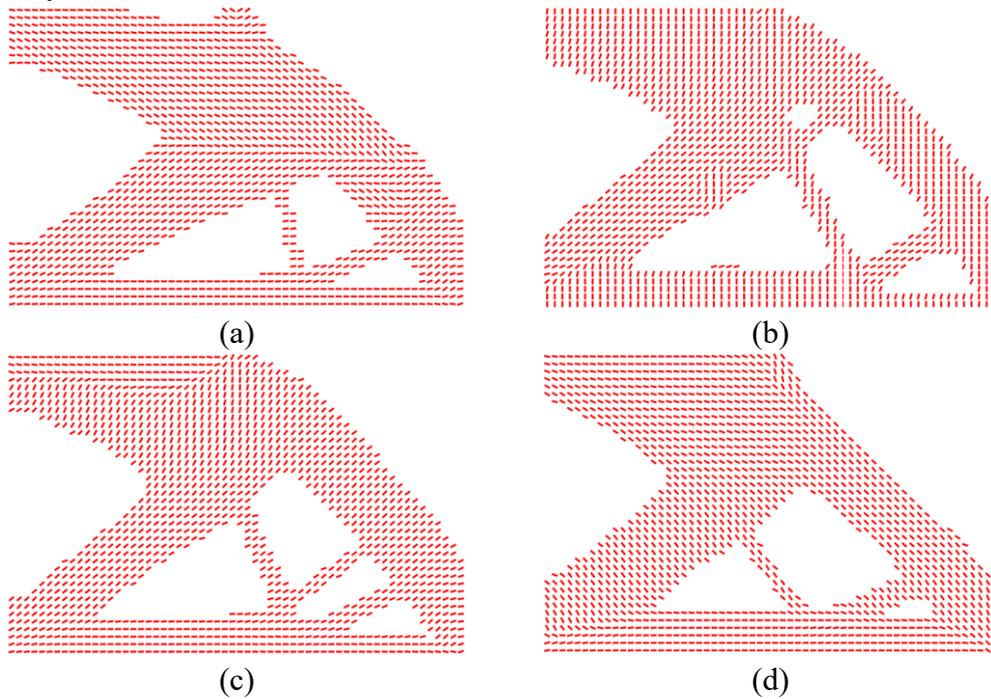

(a)          (b)

(c)          (d)

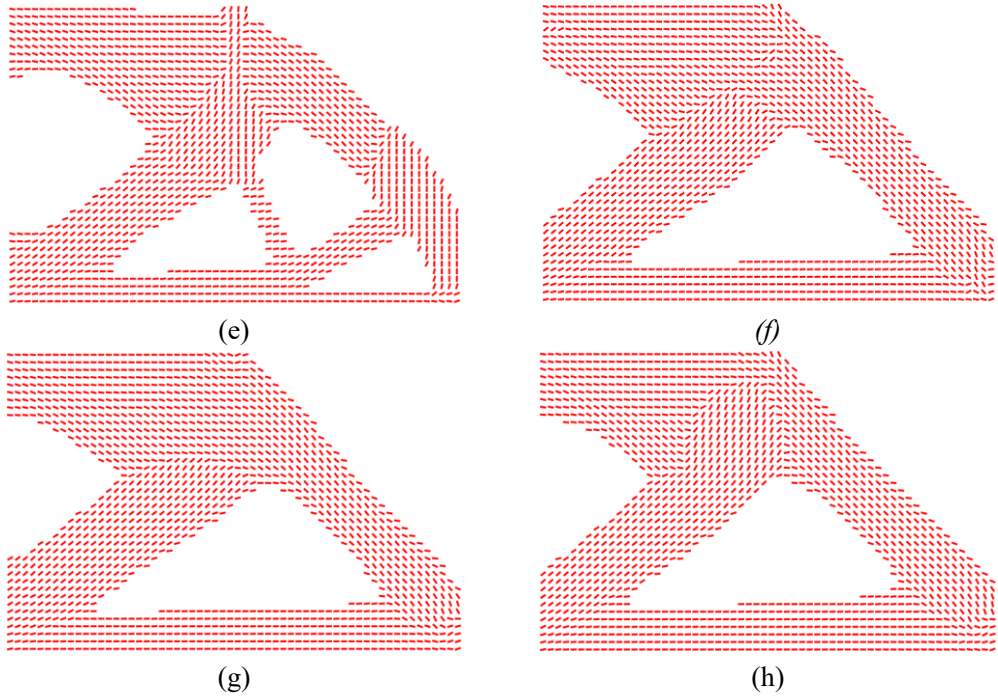

(e) (f)

(g) (h)

**Fig. 21**. The optimized fibre layout of CFAO and DSCO for Cantilever beam
with multi-load cases(a)~(h).

Fig. 22 shows the objective function and convergence rate histories of case(f) with initial fibre configuration $\theta^e_{initial}=[0°,-45°,45°,90°]$ by using the DSCO for Cantilever beam. In the optimized structure figures, areas with different colours represent different fibre angles. The total iteration is 255 and the optimized compliance reaches a low value 36.91. The first step, DMO, stops when the convergence condition is met at the 44th iteration. The convergence rate, $h_{0.95}$, reaches 0.9325 at the end of DMO. The second step, SBPTO, lasts 121 iterations.

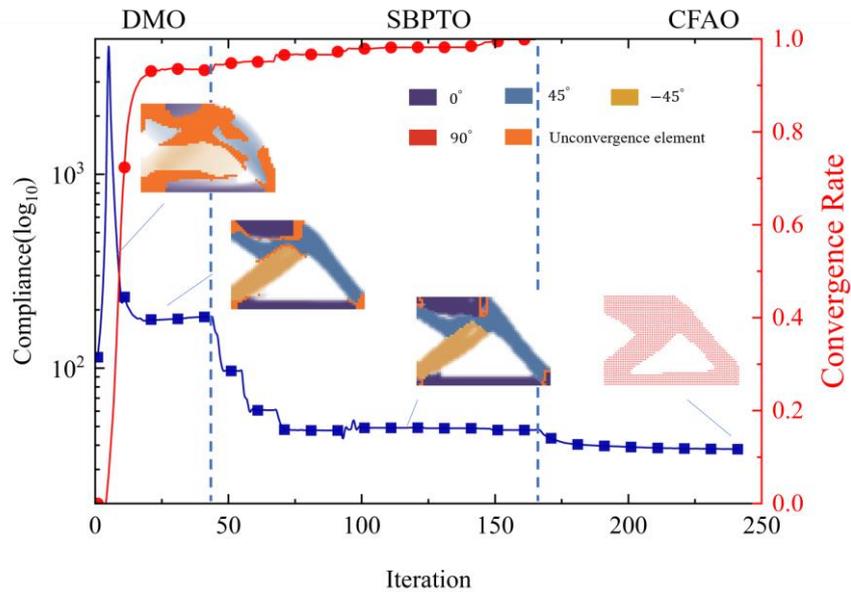

**Fig. 22**. Histories of objective function and convergence rate of case(f) for Cantilever beam with multi-load.

## 4. Conclusions

This paper establishes a new framework for concurrent optimization of topology and fibre orientation. This method aims to solve the problem of fibre convergence and make the solution closer to the optimal solution.

The main idea of this method is discrete-continuous fibre optimization. Firstly, DMO is utilized to select one fibre angle for each element from several predefined candidate angles. Secondly, since some elements cannot converge to obtain a definite angle in this process, the Sequential Binary-Phase Topology Optimization (SBPTO) is employed to improve fibre convergence rate. One specific angle orientation is treated as one material phase. The SBPTO decomposes the multiphase material topology optimization problem into a series of binary-phase material topology optimization problem. Eventually, to obtain good mechanical properties, the continuous variable of fibres orientation is designed to give more design space and use spatial filtering to make the fibre change smoothly as much as possible. The method proposed has been verified effectively by several examples. The results of numerical examples show that DSCO with $\theta_{initial}^{e} = [0°, -45°, 45°, 90°]$ has the best optimization ability in all cases. The stable optimization ability and relatively suitable computational cost make $\theta_{initial}^{e} = [0°, -45°, 45°, 90°]$ being selected as the set of predefined candidate angles.

Although the optimization framework proposed in this paper can obtain fibre angle orientations with good smoothness, the design of fibre infill pattern [21] should be introduced in order to ensure good manufacturability. Therefore, the gap between the numerical calculation results and the actual manufacturing needs to be further studied and reduced.

### Declaration of competing interest

The authors declare that they have no known competing financial interests or personal relationships that could have appeared to influence the work reported in this paper.

# Acknowledgements

We acknowledge the support provided by the Project of the National Natural Science Foundation of China (11702090) and Peacock Program for Overseas High-Level Talents Introduction of Shenzhen City (KQTD20200820113110016).